\documentclass[11pt]{amsart}
\input xypic

\usepackage{amssymb}   % For Latex2e
\usepackage{amsmath}
\usepackage{amsthm}
\usepackage{times}

\newcommand{\p}{{\mathbb P}}

\newcommand{\n}{{\mathbb N}}
\newcommand{\q}{{\mathbb Q}}
\renewcommand{\c}{{\mathbb C}}

\renewcommand{\i}{{\cal I}}

\newcommand{\ox}{{\cal O}}

\newtheorem{thm}{Theorem}

\newtheorem{cor}{Corollary}
\newtheorem{prop}{Proposition}
\theoremstyle{definition}
\newtheorem{defn}{Definition}
\newtheorem{ex}{Example}

\theoremstyle{remark}
\newcommand{\cal}{\mathcal}
\newcommand{\mult}{\operatorname{mult}}
\newcommand{\ord}{\operatorname{ord}}
\newcommand{\Exc}{\operatorname{Exc}}
\begin{document}
\title[Nef and good divisors]{A characterization of nef and good divisors by
asymptotic multiplier ideals}
\author[Francesco Russo]{Francesco Russo*}
\address{Universidade Federal de Pernambuco\\ Departamento de
Matem\' atica\phantom{aaaaaaaaaaaa}\linebreak Cidade Universitaria\\ 50670--901 Recife--PE,
Brasil} 
\thanks{*Partially  supported  by CNPq
(Centro Nacional de Pesquisa), Grant 308745/2006-0, by Edital Universal CNPq,
Grant 474475/2006-9, and  by
PRONEX/FAPERJ--Algebra Comutativa e Geometria Algebrica}
\email{frusso@dmat.ufpe.br}

\thanks{{\it AMS  Mathematics Subject Classification}:  14C20,
14F05}
\thanks{ {\it Keywords}: nef and good divisor, asymptotic
multiplier ideal, almost base point free divisor.}

\begin{abstract}
A characterization of nef and good divisors is given:   a divisor
$D$ on a smooth  complex projective  variety is nef and good if and
only if the asymptotic multiplier ideals of sufficiently high
multiples of $e(D)\cdot D$ are trivial, where $e(D)$ denotes the
exponent of the divisor $D$. Some results of the same kind are proved in the analytic setting.
\end{abstract}
\maketitle

\bigskip
\section*{Introduction}

Let $X$ be a smooth complex quasi-projective variety. One can
associate to a $\q$-divisor $D$ its multiplier ideal sheaf ${\cal
I}(D)\subseteq\ox_X$ whose zero set is the locus at which the
pair $(X,D)$ fails to have log-terminal singularities, see \cite[II.9]{Laz} and \S 1 for definitions and notation. The
multiplier ideal  ${\cal I}(D)$ reveals how bad are the
singularities of $D$. To reflect properties of the base locus of
the linear systems $|nD|$ for $n$ sufficiently large
the notion of asymptotic multiplier ideal has been introduced: the smaller
the asymptotic multiplier ideals, the worse the  asymptotic base
locus of $D$, see \cite[II.11]{Laz} for definitions and also \S 2. These two concepts and their analytic analogues,
which originated the whole theory, play an important role in
"correcting" some line bundle in order to have vanishing of
cohomology. One can consult \cite{Dem2} and \cite[II.9, II.10, II.11]{Laz} for
many applications of the theory of multiplier ideals
in analytic and algebraic geometry including results of  Lelong,
Skoda,  Siu, Nadel, Demailly, Ein and Lazarsfeld, in chronological
order, and also for complete lists of references.

In \cite[II.11.2.18]{Laz} it was shown that for a big divisor $D$ on a
smooth complex projective variety $X$ nefness is equivalent to
the triviality of the asymptotic multiplier ideals of the linear
series $|nD|$ for $n$ sufficiently large. The proof  in {\it loc. cit} is obtained
via Nadel's Vanishing Theorem for asymptotic multiplier
ideals, global generation of asymptotic multiplier ideals, \cite[II.11.2.13]{Laz}, and
boundedness of multiplicities of base loci of nef and big divisors.

Here we prove that if $D$ is a divisor on a smooth projective
complex  variety $X$ such that $\kappa(X,D)\geq 0$ and if $e(D)$
denotes the exponent of $D$, then $D$ is nef and good if and only
if the asymptotic multiplier ideals of sufficiently high
multiples of $|e(D)\cdot D|$ are trivial, i.e. if and only if
${\cal I}(n||e(D)\cdot D||)={\cal O}_X$ for $n$ sufficiently
large, Theorem \ref{char} (see also \cite[II.11.2.20]{Laz}). This generalization shows that the above condition
captures the nefness of $D$ and a sort of boundedness of the
multiplicities of the fixed components of $|nD|$ as $n$ goes to
infinity.

In the last section we recall the analytic definitions of multiplier ideal sheaf and
analytic asymptotic multiplier ideal sheaf. After analyzing the
relations between the algebraic and analytic settings, we show
by an example that the triviality of the analytic asymptotic
multiplier ideal implies nefness but not necessarily goodness.
Thus there does not exist an analytic characterization of nefness and
goodness by analytic asymptotic multiplier ideals because of the
existence of "virtual" sections, see \cite[\S 1]{DEL}.

\section{Notation and definitions}

Let $X$ be a normal complex projective variety and let $D$ be a
Cartier divisor on $X$. In \cite{G} Goodman introduced the
following definition.
\medskip

\begin{defn}(Almost base point free divisor).
A divisor $D$ is said to be {\it almost base point free} if $\forall \epsilon>0$
and $\forall  x\in X$ (not necessarily closed) there exists $n=n(\epsilon,x)$
and $D_n\in |nD|$ such that $\mult_x(D_n)<n\epsilon$.
\end{defn}
\medskip

\begin{defn}
A divisor $D$ is said to be $nef$ if $(D\cdot C)\geq 0$ for every
irreducible curve $C\subset X$.
\end{defn}
\medskip

For a nef divisor $D$ we can define the {\it numerical dimension
of $D$}: $$\nu(X,D):=\sup\{\nu\in\n: D^{\nu}\not\equiv 0\}.$$
\medskip

It is not difficult to see that $\dim(X)\geq\nu(X,D)\geq\kappa(X,D)$, where
$\kappa(X,D)$ is the Kodaira dimension of $D$.
\medskip

\begin{defn}
A nef divisor $D$ is said to be {\it good} if $\nu(X,D)=\kappa(X,D)$.
An arbitrary divisor is said to be {\it big} if $\kappa(X,D)=\dim(X)$.
\end{defn}
\medskip

By the above inequality,   a nef and big divisor is good. Let us describe some 
examples to clarify the above definitions and to put in evidence some of the
significant properties of nef and good divisors.
\medskip

\begin{ex}(A nef but not good divisor).
Let $D$ be an irreducible curve on a smooth projective surface
$S$ such that $D^2=0$ and $(D\cdot C)>0$ for every
irreducible curve $C\subset S$ with $C\neq D$. Then $D$ is a nef
divisor with $\nu(S,D)=1$ and such that $|nD|=nD$ for every
$n\geq 1$, i.e. $\kappa(S,D)=0$.

To construct explicit examples one can take as $S$ the blow-up
of $\p^2$ in $d^2$ points, $d\geq 3$, which are general on
a smooth curve $H\subset \p^2$ of degree $d$ and take as
$D$ the strict transform of $H$.

Another well known example is constructed by taking
as $D$ the "zero section" of $S=\p_F({\cal E})$, where
${\cal E}$ is the a rank 2 locally free sheaf
on an elliptic curve $F$, corresponding to (the unique)
non-splitting extension $0\to{\cal O}_F\to{\cal E}\to{\cal O}_F\to 0$.
\end{ex}
\medskip

In \cite{G} it is shown that an almost base point free divisor is
nef, see also \cite[II.11.2.19]{Laz}. The connection between the above definitions is given by the
following result which is a consequence of a theorem proved by
Kawamata in \cite[\S 2]{Kaw} (see also \cite{Mor} and
\cite{MR}).
\medskip

\begin{thm}
Let $D$ be a Cartier divisor on a complete normal complex variety
$X$. Then $D$ is almost base point free if and only if $D$ is nef
and good.
\end{thm}
\medskip

Let us recall the definitions of multiplier ideal sheaf associated
to an effective $\q$-divisor $D$ on a smooth complex projective variety $X$,
see also \cite[II.9]{Laz}. Let
$\mu:X'\to X$
be a log-resolution of $D$ and let $\Exc(\mu)$ be the sum of the exceptional divisors  of
$\mu:X'\to X$. For a $\q$-divisor $D=\sum \alpha_i D_i$ with $\alpha_i\in\q$, we
denote by $[D]=\sum[\alpha_i]D_i$ {\it the integral part of $D$}, where $[\alpha_i]$
is the integral part of $\alpha_i\in \q$.
\medskip

\begin{defn}
The {\it multiplier ideal sheaf} $${\cal I}(D)\subseteq{\cal
O}_X$$ associated to $D$ is defined to be $${\cal
I}(D)=\mu_*({\cal O}_{X'}(K_{X'/ X}-[\mu^*(D)])),$$
where $K_{X'/X}=K_{X'}-\mu^*K_X$ is the relative canonical divisor.
\end{defn}
\medskip

The multiplier ideal sheaf of $D$ does not depend on the log-resolution 
of $D$, see for example \cite[II.9.2.18]{Laz}.  Let now $|V|\subseteq|D|$ be a linear system on $X$ and let
$\mu:X'\to X$ be a log-resolution of $|V|$, i.e.
$\mu^*(|V|)=|W|+F$ where $F+\Exc(\mu)$ is a divisor with simple
normal crossing support and $|W|$ is a base point free linear
system, \cite[II.9.1.11]{Laz}.
\medskip

\begin{defn}
Fix a positive rational number $c>0$. The {\it multiplier ideal
corresponding to $c$ and $|V|$} is $${\cal I}(c\cdot
|V|)=\mu_*({\cal O}_{X'}(K_{X'/X}-[c\cdot F])).$$
\end{defn}
\medskip

Let $D$ an integral Cartier divisor on $X$ with $\kappa(X,D)\geq
0$ and let $e(D)$ be {\it the exponent of $D$}, which is by definition
the g.c.d. of the semigroup of integers $N(D)=\{ m\geq 0\,:\, |mD|\neq\emptyset\}.$ Thus there
exists a least  integer $n_0(D)$, {\it the Iitaka threshold of $D$}, such
that for every $n\geq n_0(D)$ with $e(D)| n$, $|n  D|\neq \emptyset$, see also \cite[II.11.1.A]{Laz}.
\medskip

\begin{defn}
The {\it asymptotic multiplier ideal sheaf associated to $c$ and $|D|$},
$${\cal I}(c\cdot ||D||)\subseteq {\cal O}_X,$$
is defined to be the unique maximal member among the family
of ideals $$\{{\cal I}(\frac{c}{p}\cdot |p\cdot e(D)\cdot D|)\}_{p\cdot e(D)\geq n_0(D)}.$$
\end{defn}
\medskip

In \cite[II.11.1.A]{Laz} it is shown that there exists a maximal member in the
above family, that it is unique and also that ${\cal I}(n||D||)={\cal I}(||nD||)$.
\medskip

In the next section we need  the following properties of multiplier
ideals. Let us remember that if $D=\sum \alpha_i\cdot D_i$ is a $\q$-Cartier divisor
and that if $x\in X$, then $\mult_x(D):=\sum \alpha_i\cdot \mult_x(D_i)$.
\medskip

\begin{prop}\label{Skoda}
Let $D$ be an effective $\q$-divisor on $X$. Suppose there exists a point
$x\in X$ such that $\mult_x(D)<1$. Then the multiplier ideal ${\cal I}(D)$
is trivial at $x$, i.e. ${\cal I}(D)_x={\cal O}_{X,x}$. 

If there exists
a point $x\in X$ such that $\mult_x(D)\geq\dim(X)+n-1$ for some integer $n\geq 1$, then ${\cal
I}(D)_x\subseteq m_{X,x}^n\subset {\cal O}_{X,x}$
\end{prop}
\medskip

For a proof of the first part see \cite{EV} or \cite[II.9.5.13]{Laz}.  The last 
part is proved in \cite[II.9.3.2]{Laz}. These  are algebraic
versions of  analytic results of Skoda, see \cite{Sk} or \cite[Lemma 5.6]{Dem2}.

\section{Characterization of nefness and goodness by Asymptotic Multiplier Ideals}

\begin{thm}\label{char}
Let $D$ be a  divisor on a smooth proper complex variety $X$
such that $\kappa(X,D)\geq 0$ and let $e(D)$ be the exponent
of $D$. Then $D$ is nef and good if and only if
${\cal I}(n|| e(D) D||)={\cal O}_X$ for $n$ sufficiently
large.
\end{thm}

\begin{proof} By replacing $D$ with $e(D)\cdot D$, we can assume  $e(D)=1$. 
Let us assume that $D$ is not nef and good. By Theorem 1 there exist
$\epsilon>0$ and $x\in X$, which we can assume to be a closed
point, such that for every $m\geq 1$ and for every $D_m\in|mD|$
we have $\mult_x(D_m)\geq m\epsilon$. Choose $n$ such that
$[n\epsilon]\geq \dim(X)$ and let $k$ be a sufficiently large
integer such that ${\cal I}(|| nD||)= {\cal
I}(\frac{1}{k}|knD|)$. Let $\mu:X'\to X$ be a log-resolution of
$|knD|$ constructed by first  blowing-up $X$ at $x$. The exceptional divisor
of this blow-up determines a prime divisor  $E\subset X'$ such that 
 $\mult_x(D_{kn})=\ord_E(\mu^*(D_{kn}))$
for every $D_{kn}\in|knD|$ and $\ord_E(K_{X'/ X})=\dim(X)-1$. By
definition we have $\mu^*(|knD|)=|W|+F_{kn}$ with $|W|$ base
point free. Therefore  $ kn\epsilon\leq \ord_E(F_{kn})$ and
$$\ord_E(K_{X'/ X}-[\frac{1}{k}F_{kn}]) \leq
\dim(X)-1-[n\epsilon]\leq -1,$$ yielding ${\cal I}(||nD|| )_x=
\mu_*({\cal O}_{X'}(K_{X'/X}-[\frac{1}{k}F_{kn}]))_x\subseteq \mu_*({\cal O}_{X'}(-E))_x=m_{X,x}$.
This proves that ${\cal I}(n|| e(D) D||)={\cal O}_X$ for $n$ sufficiently large implies that
$D$ is nef and good.

To prove the other implication, 
let us assume that there exists a point $x\in X$ such that ${\cal
I}(||nD||)\subseteq m_{X,x}$ for some $n\geq 1$. For $k$ sufficiently
large we have that ${\cal I}(||nD||)={\cal I}(\frac{1}{k} |knD|)$
and let $D_{kn}\in|knD|$ be a general divisor. It follows
from Proposition \ref{Skoda} that $\mult_x(D_{kn})\geq k$. Thus $D$ is not almost base point free and hence not nef and good by 
Theorem 1.
\end{proof}
\medskip

We remark that if $D$ is semiample, i.e. some
multiple of $D$ is base point free,  then $D$ is nef and good so that ${\cal I}(c\cdot ||nD||)$
is trivial for every $n\geq 1$ and for every rational number $c>0$. Moreover if $D$ is semiample,
the associated
graded algebra $R(X,D)=\oplus H^0(X,{\cal O}_X(nD))$ is finitely generated
over the base field by a result of Zariski, see \cite{Zariski} and  also \cite[I.2.1.30]{Laz}.

A nef and good divisor on a complete
normal variety is semiample if and only if the associated
graded algebra $R(X,D)$ is finitely generated
over the base field, see for example \cite{MR}
(and also  \cite{Zariski}, \cite{Wilson}, \cite{R},  \cite[I.2.3.15]{Laz} for the case of  nef and big divisors). 
Hence the triviality of  the asymptotic multiplier
ideals of  sufficiently high  multiples of a divisor
controls the nefness of the
divisor and a sort of  boundedness
of the
fixed components but not semiampleness.
In \cite{MR} and \cite{R} one can find well known
examples of nef and good divisors $D$ for which every multiple
$|nD|$ has base locus.

\section{Analytic analogue}
Let now $X$ be a compact complex manifold and let $\cal L$ be a
line bundle on $X$. Let us recall some definitions. The notation is the same as in \cite{Dem2}.
\medskip

\begin{defn} Let $\phi$ be a plurisubharmonic function, briefly a psh function, on an open subset
$\Omega\subseteq X$. The {\it Lelong mumber of $\phi$ in
$x\in\Omega$} (or of the hermitian metric $h$ having local
expression $e^{-2\phi}$) is $$\mu(\phi,x)=\liminf_{z\to
x}\frac{\phi(z)}{\log(|z-x|)}.$$
\end{defn}
\medskip

For a singular metric $h=e^{-2\phi}$ on $\cal L$ associated
to an effective divisor $D\in|\cal L|$ we have $\mu(\phi,x)=
\mult_x(D)$, that is  the  Lelong number is the analytic analogue
of multiplicity.
\medskip

The algebraic case suggests the following
definition.
\medskip

\begin{defn}
A line bundle $\cal L$ on $X$ is said {\it analytic almost base
point free} if $\forall \epsilon > 0$ and $\forall x\in X$ there
exists a possibly singular hermitian metric $h=e^{-2\phi}$  on
$\cal L$,  positive in the sense of currents (that is $\frac{i}{2\pi}\Theta_h=dd^c\phi\geq 0$ as a current) and
for which $\mu(\phi,x)<\epsilon$.
\end{defn}
\medskip

In the analytic case we have the following definition of nefness.
\medskip

\begin{defn} {\rm (\cite{DPS})}
Let $\cal L$ be a line bundle on a complex compact 
manifold $(X,\omega)$, where $\omega$ is a hermitian metric on $X$. Then $\cal L$ is said to be {\it nef} if $\forall \epsilon
> 0$ there exists a smooth hermitian metric $h_{\epsilon}$ on
$\cal L$ such that $i\Theta_{h_{\epsilon}}\geq -\epsilon\omega$.
\end{defn}

It is easy to see that if $X$ is projective, then the above
definition of nefness is equivalent to the previous one, see for example
\cite[Proposition 6.2]{Dem2}. The above condition does not imply the existence of a
smooth metric with non negative curvature on $\cal L$, see
\cite[Example 1.7]{DPS}.

If $(X,\omega)$ is a compact complex K\" ahler manifold, Demailly
proved in  \cite{Dem1} the following result which is a
generalization of Proposition 8 in \cite{G} recalled above in the algebraic setting.
\medskip

\begin{prop}{\rm (\cite{Dem1})} Let $\cal L$ be an analytic
almost base point free line bundle
on the compact complex K\" ahler manifold $(X,\omega)$. Then
$\cal L$ is nef.
\end{prop}
\medskip

On a compact K\" ahler manifold $(X,\omega)$ one defines, exactly  as in
the algebraic case, the notions of Kodaira dimension,
$\kappa(X,\cal L)$, of a line bundle $\cal L$  and, for the nef ones, of
  numerical dimension, $\nu(X,\cal L)$, see \cite[\S 6]{Dem2}. Then  as in the
algebraic case,  we have $\dim(X)\geq\nu(X,{\cal L})\geq \kappa(X,{\cal L})$
(see \cite{Dem2}), and one says that a nef line bundle is {\it good}
if $\nu(X,{\cal L})=\kappa(X,\cal L)$.
\medskip

We have the following relation between the notions of almost base
point freeness and analytic almost base point freeness, which is a
consequence of the definitions.

\begin{prop}
Let $X$ be a compact complex projective manifold and
let $\cal L$ be an almost base point free line bundle on $X$.
Then $\cal L$ is analytic almost base point free.
\end{prop}
\medskip

We now show that the converse does not hold by recalling an example
of \cite{DEL}.

\begin{ex}(An analytic almost base point free line bundle is not
almost base point free). Let $\cal E$ be a unitary flat vector
bundle on a smooth projective variety $Y$ such that no non-trivial
symmetric power of $\cal E$ or $\cal E^*$
 has sections (such vector bundles exist if for example  $Y$ is a curve
of genus $\geq 1$) and set ${\cal F}={\cal O}_C\oplus \cal E$.
Take now $X=\p(\cal F)$ and ${\cal L}={\cal O}_{\p_Y({\cal
F})}(1)$. Then for every $m\geq 1$, ${\cal L}^m$ has a unique non
trivial section which vanishes to order $m$ along the {\it
divisor at infinity} $H\subset \p(\cal F)$. Then $\cal L$ is a
nef line bundle which has a smooth semipositive metric induced by
the flat metric on $\cal E$, so that it is  analytic almost base
point free but clearly not almost base point free.
\end{ex}
\medskip

Let us remark that in  Example 2 we have  ${\cal
I}(||{\cal L}^m||)=\ox_X(-mH)$
for every $m\geq 1$ so that  analytic base point freeness cannot be characterized by
the vanishing
 of the (algebraic) asymptotic multiplier
ideal of sufficiently high multiples of $\cal L$.
\medskip

This example suggests that there should  exist an analytic
analogue of the notion of asymptotic multiplier ideal reflecting
the "boundedness of the singularities of the hermitian metrics"
on  multiples of $\cal L$. To prove this result in Proposition \ref{Skodanalitico}
we introduce the definitions of   analytic
multiplier ideal and  of  {\it metric with minimal singularities} of a line bundle on a compact
complex manifold, following Demailly, \cite{Dem2}.

\begin{defn}
Let $\phi$ be a psh function on an open subset $\Omega\subset
X$ of a complex manifold X. We associate to $\phi$ the ideal
subsheaf ${\cal I}(\phi)\subset \ox_{\Omega}$ of germs of
holomorphic functions $f\in\ox_{\Omega,x}$ such that
$|f|^2e^{-2\phi}$ is integrable with respect to the Lebesgue
measure in some local coordinates around $x$. Then ${\cal
I}(\phi)$ is said to be the {\it analytic multiplier ideal sheaf} of $\phi$.
\end{defn}
\medskip

In the analytic setting  for a line bundle $\cal L$ on $X$ whose first Chern class
lies in the closure of the cone of effective divisors, i.e.
for  a {\it pseudoeffective} line bundle,  
the notion of {\it singular hermitian metric with minimal singularities}
$h_{\min}$ can be defined in the following way (for more details see \cite[\S 13]{Dem2}
and especially the proof of \cite[Theorem 13.1.2]{Dem2}).

\begin{defn} Let $\cal L$ be a pseudoeffective  line bundle on $X$, let $h_{\infty}$
be any smooth hermitian metric on $\cal L$ and let $u=i\Theta_{h_{\infty}}(\cal L)$. Then
$$h_{\min}=h_{\infty}e^{-\psi_{\max}},$$
where $$\psi_{\max}(x)=\sup\{\psi(x)\; :\psi\; usc,\; \psi\leq 0,\;
i\partial\overline{\partial}\log (\psi)+u\geq 0\}.$$
Then one defines the {\it analytic asymptotic multiplier ideal sheaf}
of $\cal L$ as the analytic multiplier ideal sheaf of $h_{\min}$,
which will be indicated  by ${\cal I}(h_{\min})$.
\end{defn}
\medskip

The following result  follows from the fact that we have ${\cal
I}(h_{\min})=\ox_X$ if and only if $\cal L$ is analytic almost
base point free by the same argument used in the proof of Theorem
2.  We simply replace Proposition 1 by the analytic analogue for
Lelong numbers proved by Skoda,  see \cite{Sk} or \cite[Lemma 5.6]{Dem2}.

\begin{prop}\label{Skodanalitico}
Let $\cal L$ be a line bundle on a compact complex manifold $X$.
Then $\cal L$ is analytic almost base point free if and only
if $\i(h_{\min})=\ox_X$ if and only if for every point
of $X$ the Lelong numbers of $h_{\min}$ are zero.
\end{prop}

The following is a generalization 
of \cite[Proposition 13.1.4]{Dem2}.

\begin{cor}Let $\cal L$ be a nef and good line bundle on
a compact complex projective manifold  $X$. Then for $m$ sufficiently
 large $\i(||{\cal L}^{e({\cal L})m}||)={\cal I}
(h_{\min})=\ox_X$.
\end{cor}

Let us remark that  on a compact complex K\"
ahler manifold $X$ the condition  ${\cal I}(h_{\min})=\ox_X$
implies nefness but not necessarily goodness, i.e. it
does not exist an analytic characterization of nefness and
goodness by analytic asymptotic multiplier ideals.

By the above results we know that for a nef and good  line bundle on
a compact complex projective manifold the 
algebraic asymptotic multiplier ideal of its high multiples  and
its analytic asymptotic multiplier ideal  are both  trivial so that
they coincide. 

For arbitrary sections $s_1,\ldots,s_N\in H^0(X,{\cal L}^m)$ we
can take as an admissible $\psi$ function,
$\psi(x)=\frac{1}{m}\log\sum_j||\sigma_j(x)||^2_{h_{\infty}}+C$, with $C$ a costant.
From this it follows, see \cite{DEL} and \cite{Dem2}, that if $X$
is a complex compact projective manifold and if $\kappa(X,{\cal
L})\geq 0$, then ${\cal I}(||{\cal L}^m||)\subseteq {\cal
I}(h_{\min}).$
Example 2 above shows that the inclusion
$\i(||{\cal L}^m||)\subseteq\i(h_{\min})$ can be strict.

One conjectures that for arbitrary big line bundles the asymptotic algebraic multiplier
of its multiples and its analytic multiplier ideal coincide, see \cite{DEL} and also \cite[II.11.1.11]{Laz}.

\section*{Acknowledgements}

This note originated several years ago  participating to the  {\it School on Vanishing Theorems and Effective Results
in Algebraic Geometry}, ICTP Trieste, 2000, supported by a grant of ICTP. With a long delay, I would like  to thank 
the organizers for inviting me  and  
all the lecturers for the stimulating courses and for the very interesting  lecture notes.

Last but not least, I am  grateful to the organizers of the Ghent Conference {\it Linear Systems and Subschemes},
april 2007, for inviting me and for the nice and stimulating atmosphere during the meeting
 and finally to Ghent University for financial support.

\end{document}